\documentclass[10pt]{article}

\usepackage{theorem,amssymb,amsmath,url}

\usepackage{graphicx}

\usepackage{hyperref}

\usepackage{color}                    
\definecolor{red}{rgb}{1,0,.2}        
\definecolor{cjp}{rgb}{.1,.7,.2}        
\definecolor{fmdc}{rgb}{1,0,.8}        
\definecolor{gree}{rgb}{0,0.7,0.2}        
\definecolor{blue}{rgb}{0,0,.9}

\newcommand{\morf}[2]{\Big\{\begin{aligned} 0 & \rightarrow #1\\[-.1cm] 1 & \rightarrow #2  \end{aligned}}

\newcommand{\T}{x_{3/2}}

\newtheorem{theorem}{Theorem}

\newtheorem{proposition}[theorem]{Proposition}
\theorembodyfont{\rm}

\advance \topmargin by -\headheight
\advance \topmargin by -\headsep
\evensidemargin \oddsidemargin
\begin{document}

{\huge \bf Two-block substitutions and morphic words}

\begin{center}

{\small Version 2, February 20, 2023}


\bigskip

{\Large \bf  Michel Dekking and Michael Keane}

\medskip

{\small f.m.dekking@tudelft.nl, Michel.Dekking@cwi.nl}\\
{\small m.s.keane@tudelft.nl}

\end{center}





\begin{abstract}
 \noindent We consider in general two-block substitutions and their fixed points. We prove that some of them have a simple structure: their fixed points are morphic sequences. Others are intrinsically more complex, such as the Kolakoski sequence. We prove this for the Thue-Morse sequence in base 3/2.
\end{abstract}

\medskip



\bigskip

\section{Introduction}

 \qquad Let $A=\{0,1\}$, $A^*$ the monoid of all words over $A$, and let  $T^*$ be the submonoid of 0-1-words of even length.
 A \emph{two-block substitution} $\kappa$ is a map
 $$\kappa:\: \{00,01,10,11\}\rightarrow A^*.$$
 A two-block substitution $\kappa$ acts on $T^*$ by defining for $w_1w_2\dots w_{2m-1}w_{2m}\in T^*$
 $$\kappa(w_1w_2\dots w_{2m-1}w_{2m}) = \kappa(w_1w_2)\dots \kappa(w_{2m-1}w_{2m}).$$
 In the case that $\kappa(T^*)\subseteq T^*$, we call $\kappa$ \emph{$2$-block stable}. This property entails that the iterates $\kappa^n$ are all well-defined for $n=1,2,\dots$.

 The most interesting example of a two-block substitution that is \emph{not} 2-block stable is the Oldenburger-Kolakoski two-block substitution $\kappa_{\rm\sc K}$ given by
 $$\kappa_{\rm\sc K}(00)=10,\quad \kappa_{\rm\sc K}(01)=100, \quad\kappa_{\rm\sc K}(10)=110,\quad \kappa_{\rm\sc K}(11)=1100.$$

 The fact that $\kappa_{\rm\sc K}$ is not 2-block stable, and so its iterates $\kappa_{\rm\sc K}^n$ are not defined, makes it very hard to establish properties of the fixed point $x_{\rm\sc K}=110010\dots$  (usually written as 221121...)  of $\kappa_{\rm\sc K}$, see, e.g., \cite{Dekk-Kol}.

 \medskip

In Section \ref{sec:gold} we show that even if a two-block substitution  $\kappa_{\rm\sc K}$ is {\it not} 2-block stable, then still it can be well-behaved in the sense that its fixed points are pure morphic words.


In Section \ref{sec:TMK} we prove that the Thue-Morse word in base 3/2 is not well-behaved: it can not be generated as a coding of a fixed point of a morphism.

\section{Two-block substitutions with conjugated morphisms}\label{sec:gold}

 Let $\kappa$ be a two-block substitution on $T^*$, and let  $\sigma$ be a morphism on $A^*$ with $\sigma(T^*)\subseteq T^*$. We say $\kappa$ and $\sigma$ commute if  $\kappa\sigma(w)=\sigma\kappa(w)$ for all $w$ from $T^*$.\\
 In this case we say that $\sigma$ is {\it conjugated} to $\kappa$.\\
Note that if  $\kappa\sigma=\sigma\kappa$,  then for all $n\ge1$ one has $\kappa\sigma^n=\sigma^{n}\kappa$ on  $T^*$.

\medskip

Let $\sigma: A^*\rightarrow A^*$ be a morphism. Then $\sigma$ induces a two-block substitution $\kappa_\sigma$ by defining
 $$\kappa_\sigma(ij)=\sigma(ij)\quad {\rm for\:} i,j\in A.$$

 We mention the following property of $\kappa_\sigma$, which is easily proved by induction.

 \begin{proposition}\label{prop:iter}
 Let $\sigma: A^*\rightarrow A^*$ be a morphism, let $n$ be a positive integer, and suppose that $\kappa_\sigma$ is two-block stable.
 Then $\kappa_\sigma^n=\kappa_{\sigma^n}$.
 \end{proposition}

 We call $\sigma$ the  {\it trivial conjugated morphism} of the two block substitution $\kappa_\sigma$.

 Not all morphisms $\sigma$ can occur as trivial conjugated morphisms, but many will be according to the following simple property.

\begin{proposition}\label{prop:conj}
 Any morphism $V$ on $\{0,1\}$ with the lengths of $\sigma(0)$ and $\sigma(1)$ both odd or both even is conjugated to the two-block substitution $\kappa=\kappa_\sigma$.
 \end{proposition}

 Example:  for the Fibonacci morphism $\varphi$ defined by $\varphi(0)=01, \varphi(1)=0$, one can take the third power $\varphi^3$ to achieve this (cf.~\cite[A143667]{OEIS}). 

 \medskip

 In the remaining part of this section we discuss non-trivial conjugated morphisms.

\begin{theorem}\label{th:commu} Let $\kappa$ be a two-block substitution on $T^*$ conjugated with a morphism $\sigma$ on $A^*$. Suppose that there exist $i,j$ from $A$ such that  $\kappa(ij)$  has prefix $ij$, and such that $ij$ is also prefix of a fixed point $x$ of $\sigma$. Then also $\kappa$ has fixed point $x$.
\end{theorem}

\medskip

\noindent \emph{Proof:} Letting $n\rightarrow \infty$ in $\kappa\sigma^n(ij)=\sigma^n\kappa(ij)=\sigma^n(ij\dots)$ gives $\kappa(x)=x$.\qquad$\Box$

\bigskip

The Pell word $w_{\rm \sc P} = 0010010001001\dots$ is the unique fixed point of the Pell morphism $\pi$ given by
 $$\pi: \; \morf{001}{0}.$$  

 The following result proves a conjecture from R.J.~Mathar in \cite[A289001]{OEIS}. The difficulty here is that since the 2-block substitution in Theorem \ref{th:Pell} has the property that $\kappa(0010)=0010010$ has odd length, the two-block substitution  $\kappa$ is not 2-block stable.

\begin{theorem}\label{th:Pell}  Let $\kappa$ be the two-block
substitution\footnote{Here it is not necessary to define $\kappa(11)$, since 11 does not occur in images of words without 11 under $\kappa$.}:\\[-.3cm]
$$\kappa: \; \Bigg\{\begin{aligned} 00 & \rightarrow 0010\\[-.1cm] 01 & \rightarrow 001  \\[-.1cm] 10 & \rightarrow 010. \end{aligned}$$
Then the  unique fixed point of $\kappa$ is the Pell word $w_{\rm \sc P}$.
\end{theorem}

\noindent \emph{Proof:} We apply Theorem \ref{th:commu} with $ij=00$.\\
Note first that $\pi(T^*)\subseteq T^*$. Next, we have to establish that  $\kappa$ and $\pi$ commute on  $T^*$.

It suffices  to check this for the three generators 00, 01 and 10 from the four generators of  $T^*$:\\[-.6cm]
\begin{eqnarray*}
\kappa\pi(00)&=& \kappa(001001) = 0010010001 = \pi(0010)=\pi\kappa(00),\\
\kappa\pi(01)&=&\kappa(0010) = 0010010 = \pi(001)=\pi\kappa(01),\\
\kappa\pi(10)&=& \kappa(0001) = 0010001 = \pi(010)=\pi\kappa(10).\qquad\Box
\end{eqnarray*}

\section{Thue-Morse in base 3/2}\label{sec:TMK}

 A natural number $N$ is written in base 3/2 if $N$ has the form
 \begin{equation}\label{eq:exp}
    N= \sum_{i\ge 0} d_i\Big( \frac32\Big)^i,
 \end{equation}  \vspace*{-.0cm}
  with digits $d_i=0,1$ or 2.

 We  write these expansions as
  $${\rm SQ}(N) = d_R(N)\dots d_1(N)d_0(N)= d_R\dots d_1d_0.$$

Let for $N\ge 0$, $s_{3/2}(N):=\sum_{i=0}^{i=R} d_i(N)$ be the sum of digits function of the base 3/2 expansions.
The Thue-Morse word in base 3/2 is the word  $(\T(N)):= (s_{3/2}(N) \!\!\! \mod 2) \,= \, 0100101011011010101\dots$

\begin{theorem}\label{th:TM} {\bf  (\cite{Dekking-3-2})} Let  the two-block substitution $\kappa_{\rm\sc TM}$ be defined by
\begin{displaymath}
\kappa_{\rm\sc TM}: \begin{cases}
00 & \rightarrow 010\\[-.1cm]
01 & \rightarrow 010  \\[-.1cm]
10 & \rightarrow 101  \\[-.1cm]
11 & \rightarrow 101
\end{cases}
\end{displaymath}
Then the  word $\T$  is the fixed point of $\kappa_{\rm\sc TM}$ starting with $0$.
\end{theorem}

\medskip

The Thue-Morse word $t$ is fixed point with prefix 0 of the Thue-Morse morphism $\tau: 0\rightarrow 01,\: 1\rightarrow 10$.
It satisfies the recurrence relations $t(2N)=t(N), \:t(2N+1)=1-t(N)$.

The fixed point $\T$  satisfies very similar recurrence relations:
$$\T(3N+1)=\T(2N+1),\; \T(3N+2)=1-\T(2N+1),\; \T(3N+3)=\T(2N+1).$$

We call $\kappa_{\rm\sc TM}$ the \emph{Thue-Morse two-block substitution}.


\medskip

We now discuss the Kolakoski word $x_{\rm\sc K}$. This word was introduced by Kolakoski (years after Oldenburger
\cite{Olden}) as a problem in \cite{Kol}. The problem was to prove that $x_{\rm\sc K}$ is not eventually periodic. Its solution in \cite{Ucoluk} is however incorrect (The claim  that words $w$ with minimal period $N$ in $www\dots$ map to words with period $N_1$ satisfying $N < N_1 < 2N$ by replacing run lengths by the runs themselves is false. For example, if the period word is $w = 21221$, then $ww$  maps to the period word 2212211211211221, or its binary complement image.)
A stronger result was proved by both Carpi \cite{Carpi} and Lepist\"o \cite{Lepisto}: $x_{\rm\sc K}$ does not contain any cubes.
The fixed point $\T$ of $\kappa_{\rm\sc TM}$ has more repetitiveness. It contains for example the fourth power 01010101.

\medskip

The Thue-Morse word is a purely morphic word, i.e., fixed point of a morphism. It is known that the Kolakoski word is not purely morphic (\cite{Dekk-Kol}). However it is still open whether the Kolakoski word is morphic, i.e., image under a coding (letter to letter map) of a fixed point of a morphism. The tool here is the subword complexity function $(p(N))$, which gives the number of words of length $N$ occurring in an infinite word. A well known result tells us that when the subword complexity function increases too fast, faster than $N^2$, then a word can not be morphic. There is one example of a two-block substitution which yields a word that is not morphic given by Lepist\"o in the paper \cite{Lepisto-comp}.

\begin{theorem}\label{th:Lep} {\bf  (\cite{Lepisto-comp})} Let  the two-block substitution $\kappa_{\rm\sc L}$ be defined by
\begin{displaymath}
\kappa_{\rm\sc L}: \begin{cases}
00 & \rightarrow 011\\[-.1cm]
01 & \rightarrow 010  \\[-.1cm]
10 & \rightarrow 001  \\[-.1cm]
11 & \rightarrow 000
\end{cases}
\end{displaymath}
Then the   fixed point $010011000011\dots$ of $\kappa_{\rm\sc L}$ has subword complexity function $p(N)$ satisfying $p(N)>C.N^t$ for some $C>0$ and $t>2$.
\end{theorem}

We do not know how the prove this 'faster than quadratic' property for the base 3/2 Thue-Morse word, but still we can use Lepist\"o's result to obtain the following.

\begin{theorem}\label{th:non-morph} The base $3/2$ Thue-Morse word $\T$ is not a morphic word.
\end{theorem}

\medskip

\noindent \emph{Proof:} The idea is to introduce the base 3/2 Toeplitz word $x_{\rm T}$ fixed point of the two-block substitution given by
\begin{displaymath}
\kappa_{\rm\sc T}: \begin{cases}
00 & \rightarrow 111\\[-.1cm]
01 & \rightarrow 110  \\[-.1cm]
10 & \rightarrow 101  \\[-.1cm]
11 & \rightarrow 100
\end{cases}
\end{displaymath}
We call this the \emph{base $3/2$ Toeplitz word} because  $x_{\rm T}(3N)=1$ for all $N$, and because it can be simply obtained from the base 3/2 Thue-Morse word---similar to the way the `classical' Toeplitz word can be obtained from the base 2 Thue-Morse  word $t$, as in \cite[Lemma 3]{Alletal}. The corresponding base 3/2 equation is as follows:
\begin{equation}\label{eq:ThT}
x_{\rm T}(N)=1-|\T(3N-2)-\T(3N+1)| \quad {\rm for}\; N=1,2,\dots.
\end{equation}
See  \cite[Theorem 3.2]{Edgar-} for a proof of Equation \ref{eq:ThT} (N.B.~There is a shift of indices). 

The crucial observation is now that the base 3/2 Toeplitz two-block substitution is just the binary complement of the $\kappa_{\rm\sc L}$ two-block substitution. In particular Theorem \ref{th:Lep} also holds for the base 3/2 Toeplitz word, and so $x_{\rm T}$ can not be a morphic word. 

 Suppose that the base 3/2 Thue-Morse word $(\T(N))$ is a morphic word. Then an application of \cite[Theorem 7.9.1]{AllShall} yields that the word $(\T(3N-2))$ is morphic. Next, \cite[Theorem 7.6.4]{AllShall} gives that the direct product word $([\T(3N-2),\T(3N+1)])$  is morphic. Finally, another application of \cite[Theorem 7.9.1]{AllShall} yields that according to Equation \ref{eq:ThT}
this direct product word maps to a morphic word $(x_{\rm T}(N))$ under the morphism $[0,0]\mapsto 1,\: [0,1]\mapsto 0,\: [1,0]\mapsto 0,\: [1,1]\mapsto 1$. But this contradicts the fact that  $(x_{\rm T}(N))$ is not morphic. Hence the  base 3/2 Thue-Morse word is not a morphic word. \qquad$\Box$


\section{Acknowledgement} We are grateful to Jean-Paul Allouche for several useful comments.

\end{document}